\documentclass[11pt]{article}
\usepackage{amsmath,amsfonts,amssymb,amsthm}
\begin{document}

\newcommand{\1}{{{\bf 1}}}
\newcommand{\id}{{\rm id}}
\newcommand{\Hom}{{\rm Hom}\,}
\newcommand{\End}{{\rm End}\,}
\newcommand{\Res}{{\rm Res}\,}
\newcommand{\Image}{{\rm Im}\,}
\newcommand{\Ind}{{\rm Ind}\,}
\newcommand{\Aut}{{\rm Aut}\,}
\newcommand{\Ker}{{\rm Ker}\,}
\newcommand{\gr}{{\rm gr}}
\newcommand{\Der}{{\rm Der}\,}

\newcommand{\Z}{\mathbb{Z}}
\newcommand{\Q}{\mathbb{Q}}
\newcommand{\C}{\mathbb{C}}
\newcommand{\N}{\mathbb{N}}
\newcommand{\g}{\mathfrak{g}}
\newcommand{\gl}{\mathfrak{gl}}
\newcommand{\h}{\mathfrak{h}}
\newcommand{\wt}{{\rm wt}\;}
\newcommand{\A}{\mathcal{A}}
\newcommand{\D}{\mathcal{D}}
\newcommand{\Lie}{\mathcal{L}}

\def \<{\langle}
\def \>{\rangle}
\def \be{\begin{equation}\label}
\def \ee{\end{equation}}
\def \bex{\begin{exa}\label}
\def \eex{\end{exa}}
\def \bl{\begin{lem}\label}
\def \el{\end{lem}}
\def \bt{\begin{thm}\label}
\def \et{\end{thm}}
\def \bp{\begin{prop}\label}
\def \ep{\end{prop}}
\def \br{\begin{rem}\label}
\def \er{\end{rem}}
\def \bc{\begin{coro}\label}
\def \ec{\end{coro}}
\def \bd{\begin{de}\label}
\def \ed{\end{de}}

\newtheorem{thm}{Theorem}[section]
\newtheorem{prop}[thm]{Proposition}
\newtheorem{coro}[thm]{Corollary}
\newtheorem{conj}[thm]{Conjecture}
\newtheorem{exa}[thm]{Example}
\newtheorem{lem}[thm]{Lemma}
\newtheorem{rem}[thm]{Remark}
\newtheorem{de}[thm]{Definition}
\newtheorem{hy}[thm]{Hypothesis}
\makeatletter \@addtoreset{equation}{section}
\def\theequation{\thesection.\arabic{equation}}
\makeatother \makeatletter

\begin{Large}
\begin{center}
\textbf{Vertex Poisson algebras associated with Courant algebroids
and their deformations; I}
\end{center}
\end{Large}
\begin{center}{
Gaywalee Yamskulna\footnote{E-mail address: gyamsku@ilstu.edu}\footnote{Partially supported by a New Faculty Initiative grant from the College of Arts and Sciences, ISU}\\
Department of Mathematical Sciences, Illinois State University, Normal, IL 61790\\
and\\
Institute of Science, Walailak University, Nakon Si Thammarat,
Thailand}
\end{center}

\begin{abstract}
This is the first of two papers on vertex Poisson algebras
associated with Courant algebroids, and their deformations. In
this work, we study relationships between vertex Poisson algebras
and Courant algebroids. For any $\N$-graded vertex Poisson algebra
$A=\coprod_{n\in\N} A_{(n)}$, we show that $A_{(1)}$ is a Courant
$A_{(0)}$-algebroid. On the other hand, for any Courant
$\mathcal{A}$-algebroid $\mathcal{B}$, we construct an $\N$-graded
vertex Poisson algebra $A=\coprod_{n\in\N}A_{(n)}$ such that
$A_{(0)}$ is $\mathcal{A}$ and the Courant $\mathcal{A}$-algebroid
$A_{(1)}$ is isomorphic to $\mathcal{B}$ as a Courant
$\mathcal{A}$-algebroid.
\end{abstract}

\section{Introduction}

A vertex Poisson algebra is an analogue of the notion of a Poisson
algebra in the category of vertex algebras. It is a combination of
a commutative vertex algebra structure (or equivalently a
differential algebra structure) and a vertex Lie algebra structure
with a natural compatibility condition (see \cite{bed, bf}). In
fact, one may consider vertex Poisson algebras as the classical
limit of vertex algebras (see \cite{bed, bf}). Important features
of the formalism of vertex Poisson algebras can be traced back to
the work of I.M. Gelfand, L.A. Dickey  and others on the
Hamiltonian structure of integrable hierarchies of solition
equations (see \cite{di}). They play a prominent role in the study
of certain connections between the classical and quantum
Drinfeld-Sokolov reductions (see \cite{bf} Chapter 15). They also
play a crucial role in the study of generating subspaces of vertex
algebras with a certain property analogous to the well known
Poincar\'{e}-Birkhoff-Witt spanning property (see \cite{li2, li3}
and references their in).
%In \cite{li2},
%motivated by classical results in Lie theory, Li introduced and
%studied a notion of what so called good increasing filtration for
%a vertex algebra $V$ and proved that the associated graded vector
%space $gr V$ of $V$ with respect to a good increasing filtration
%is naturally a vertex Poisson algebra. For any $\N$-graded vertex
%algebra $V=\coprod_{n\in\N} V_{(n)}$ with $V_{(0)}=\C{\bf 1}$, the
%vertex Poisson algebra $gr V$ was essentially used in \cite{kl,
%gn, bu1, bu2, abd, nt} in the study on generating subspaces of $V$
%with a certain property analogous to the well known
%Poincar$\acute{e}$-Birkhoff-Witt spanning property. In \cite{li3},
%for any vertex algebra $V$, Li associated a canonical decreasing
%sequence of subspaces and prove that the associated graded vector
%space $gr(V)$ is naturally a vertex Poisson algebra. He then use
%vertex Poisson algebras $gr(V)$ to study generating subspaces of
%certain types for lower truncated $\Z$-graded vertex algebras.
% In particular, by using this
%vertex Poisson algebra $gr(V)$, he proved that for any vertex
%algebra $V$, $C_2$-cofiniteness implies $C_n$-cofiniteness for all
%$n\geq 2$. The notion of $C_2$ was introduced and used in the
%fundamental study of Zhu on modular invariance (see \cite{z}).

Courant algebroids first appeared in the work of T. Courant on
Dirac structures on manifolds (see \cite{co}). They were later
studied by Liu, Weinstein, and Xu who used Courant algebroids to
generalize the notion of the Drienfeld double to Lie bialgebroids
(see \cite{liuwx}). One can think of a Courant algebroid as a
certain kind of algebroid for which the condition of
antisymmetricity, the Leibniz rule and the Jacobi identity for the
bracket are all relaxed in certain ways. Examples of Courant
algebroids are the doubles of Lie bialgebras and the bundle
$TM\oplus T^*M$. (see \cite{liux, co}). Courant algebroids play an
important role in the study of dynamical $r$-matrices (see
\cite{liux}). Also, it was shown in \cite{row} that Courant
algebroids can be viewed as strongly homotopy Lie algebras.
Furthermore, Courant algebroids have a close connection to gerbs
(see \cite{brc}) and play a significant role in the study of the
theory of vertex algebroids (see \cite{br}).

The purpose of this paper is to study relationships between
Courant algebroids and $\N$-graded vertex Poisson algebras.
Precisely, we show that one can obtain Courant algebroids from
$\N$-graded vertex Poisson algebras and vice versa. For an
$\N$-graded vertex Poisson algebra $A=\coprod_{n\in\N}A_{(n)}$,
the homogeneous subspace $A_{(0)}$ is a unital commutative
associative algebra and the homogeneous subspace $A_{(1)}$ is a
module of an associative algebra $A_{(0)}$. Moreover, the skew
symmetry and the half commutator formula of a vertex Poisson
algebra give rise to several compatibility relations. These
structures on $A_{(0)}\oplus A_{(1)}$ are summarized in the notion
of what was called a 1-truncated conformal algebra. By using skew
symmetry and the fact that a family of bilinear operations
$$_i:A\otimes_{\C}A\rightarrow A, \ \ a\otimes b\mapsto a_ib$$
are derivations of the commutative product on $A$, we show that
$A_{(1)}$ is, in fact, a Courant $A_{(0)}$-algebroid. On the other
hand, for a given Courant $\mathcal{A}$-algebroid $\mathcal{B}$,
we construct an $\N$-graded vertex Poisson algebra
$A=\coprod_{n\in\N}A_{(n)}$ with $A_{(0)}=\mathcal{A}$ and the
homogeneous subspace $A_{(1)}$ isomorphic to $\mathcal{B}$ as a
Courant algebroid. For a given 1-truncated conformal algebra
$\mathbb{A}\oplus \mathbb{B}$ over $\C$, we set
$C(\mathbb{A}\oplus \mathbb{B})= \C[\D]\otimes(\mathbb{A}\oplus
\mathbb{B})$ where $\D$ is a formal variable. We use the fact that
$\mathbb{A}\oplus \mathbb{B}$ is a 1-truncated conformal algebra
to show that a certain quotient space ${\mathcal{C}}_{\mathbb{B}}$
of $C(\mathbb{A}\oplus \mathbb{B})$ is a vertex Lie algebra. It
was shown in \cite{li2} that a symmetric algebra of a vertex Lie
algebra is a vertex Poisson algebra. Applying this result to
${\mathcal{C}}_{\mathbb{B}}$, we then obtain a vertex Poisson
algebra $S({\mathcal{C}}_{\mathbb{B}})$. In fact,
$S({\mathcal{C}}_{\mathbb{B}})$ is an $\N$-graded vertex Poisson
algebra whose degree-zero subspace is $S(\mathbb{A})$ and whose
degree-one subspace is $S(\mathbb{A})\cdot \mathbb{B}$. Here,
$S(\mathbb{A})$ is the symmetric algebra over the space
$\mathbb{A}$. For any Courant $\mathcal{A}$-algebroid
$\mathcal{B}$, $\mathcal{A}\oplus\mathcal{B}$ is a 1-truncated
conformal algebra. Hence, $S({\mathcal{C}}_{\mathcal B})$ is an
$\N$-graded vertex Poisson algebra. By using the fact that
$\mathcal{B}$ is a Courant $\mathcal{A}$-algebroid, we can show
that a certain quotient space $S_{\mathcal{B}}$ of
$S({\mathcal{C}}_{\mathcal B})$ is an $\N$-graded vertex Poisson
algebra whose the degree-zero subspace is $\mathcal{A}$ and the
degree-one subspace is isomorphic to $\mathcal{B}$ as a Courant
$\mathcal{A}$-algebroid.

This paper is organized as follows: In section 2 we recall notions
of vertex Lie algebras, and vertex Poisson algebras. In section 3,
we review the notions of 1-truncated conformal algebras and
Courant $\mathcal{A}$-algebroids, and we construct Courant
algebroids from $\N$-graded vertex Poisson algebras. In section 4,
we construct an $\N$-graded vertex Poisson algebra from any
1-truncated conformal algebra. In section 5, for any Courant
$\mathcal{A}$-algebroid $\mathcal{B}$, we construct an $\N$-graded
vertex Poisson algebra whose degree-zero subspace is $\mathcal{A}$
and degree-one subspace is isomorphic to $\mathcal{B}$ as a
Courant $\mathcal{A}$-algebroid.

We use here the standard formal variable notations and conventions
as defined in \cite{fhl,lli}. Also, we use a notation $\N$ for the
set of nonnegative integers.

\section{Vertex Lie algebras, and vertex Poisson algebras}

In this section we review the definitions of vertex Lie algebras
and vertex Poisson algebras, as well as their properties that will
be used in later sections.

Let $W$ be a vector space. Following \cite{p, li2}, for a formal
series
$$f(x_1,...,x_n)=\sum_{m_1,...,m_n\in\Z}u(m_1,...,m_n)x_1^{-m_1-1}...x_n^{-m_n-1}\in
W[[x_1^{\pm 1},...,x_n^{\pm 1}]],$$ we define
$$Sing\ \  f(x_1,...,x_n)=\sum_{m_1,...,m_n\in\N}u(m_1,...,m_n)x_1^{-m_1-1}...x_n^{-m_n-1}.$$
For $1\leq i\leq n$, we have
\begin{equation}\label{dxsing}
\frac{\partial }{\partial x_i}Sing\ \ f(x_1,...,x_n)=Sing\ \
\frac{\partial}{\partial x_i} f(x_1,...,x_n).
\end{equation}
For any nonempty subset $S=\{ i_1,...,i_k\}$ of $\{1,...,n\}$, a
formal series $f(x_1,...,x_n)$ can be viewed as a formal series
$\widetilde{f}(x_{i_1},...,x_{i_k})$ in variables
$x_{i_1},...,x_{i_k}$ with coefficients in the vector space
$V[[x_{j}^{\pm 1}|j\not\in S]]$. We define
\begin{equation}
Sing_{x_{i_1},...,x_{i_k}}\ \ f(x_1,...,x_n)=Sing\ \
\widetilde{f}(x_{i_1},...,x_{i_k}).
\end{equation}
Hence, $$Sing_{x_1}\ \
f(x_1,...,x_n)=\sum_{m\in\N,m_2,...,m_n\in\Z}u(m,m_2,...,m_n)x_1^{-m-1}x_2^{-m_2-1}...x_n^{-m_n-1}$$
and
$$Sing\ \  f(x_1,...,x_n)=Sing_{x_1}...Sing_{x_n}\ \
f(x_1,...x_n).$$

\begin{prop}\label{singsing}\cite{p} Let $W$ be a vector space and let
$B\in W((x_1,...,x_n))$, $M\in (\End W)[[x_1,...,x_n]].$ Then
$$Sing \ \ (M\cdot Sing\ \ (B))=Sing\ \ (MB).$$
\end{prop}

Next, we recall the notions of vertex Lie algebra and vertex
Poisson algebra.
\begin{de}\cite{k, p} {\em A vertex Lie algebra} is a vector space $A$
equipped with a linear operator $\partial$ and a linear map
\begin{eqnarray*}
Y_-:A&\rightarrow &\Hom(A,x^{-1}A[x^{-1}])\\
a&\mapsto&Y_-(a,x)=\sum_{n\in \N}a_nx^{-n-1}\ \ (\text{where }
a_n\in \End A)
\end{eqnarray*}
such that the following axioms hold for $a,b\in A$:
\begin{eqnarray}
& &a_nb=0\ \ \text{ for $n$ sufficiently large};\label{yhn}\\
& &Y_-(\partial a,x)=\frac{d}{dx}Y_-(a,x);\label{yhp}\\
& &Y_-(a,x)b=Sing(e^{x\partial}Y_-(b,-x)a);\label{yhs}\\
&
&{[Y_-(a,x_1),Y_-(b,x_2)]}=Sing(Y_-(Y_-(a,x_1-x_2)b,x_2)).\label{yha}
\end{eqnarray}Following \cite{p}, we call (\ref{yha}) {\em the half commutator
formula}.\end{de}
\begin{rem} In terms of components, (\ref{yhp})-(\ref{yha}) are equivalent to
\begin{eqnarray}
& &(\partial a)_nb=-na_{n-1}b,\label{hp}\\
& &a_nb=\sum_{i\geq 0}(-1)^{n+i+1}\frac{1}{i!}\partial^i
b_{n+i}a,\label{hs}\\
& &a_mb_nc-b_na_mc=\sum_{i=0}^m{m\choose
i}(a_ib)_{m+n-i}c,\label{ha}
\end{eqnarray}
for $a,b,c\in A, m,n\in \N$.
\end{rem}
\begin{rem} Any vertex algebra is a vertex Lie algebra with
$\partial=\D$. (See appendix for the definition of a vertex
algebra.)
\end{rem}
\begin{prop}\label{rhd}\cite{k} For $a\in A$,
$[\partial ,Y_-(a,x)]=\frac{d}{dx}Y_-(a,x)=Y_-(\partial a,x).$
\end{prop}

\begin{prop}\label{haf}\cite{li2} Let $A$ be a vector space equipped with a
linear operator $\partial$ and equipped with a linear map
\begin{eqnarray}
Y_-:A&\rightarrow&\Hom(A, x^{-1}A[x^{-1}])\nonumber\\
a&\mapsto&Y_-(a,x)
\end{eqnarray}
such that the following conditions hold for all $a,b\in A$:
\begin{eqnarray*}
& &[\partial, Y_-(a,x)]=\frac{d}{dx}Y_-(a,x),\\
& &Y_-(a,x)b=Sing(e^{x\partial}Y_-(b,-x)a).
\end{eqnarray*}
Then the half commutator formula for an ordered triple $(a,b,c)$
implies the half commutator formula for any permutation of
$(a,b,c)$.
\end{prop}

{\em A differential algebra} is a commutative associative algebra
$A$ with the identity 1 equipped with a derivation $\partial$. We
often denote the differential algebra by $(A,\partial )$. A subset
$U$ of $A$ generates $A$ as a differential algebra if $\partial ^n
U$ for $n\in\N$ generate $A$ as an algebra.
\begin{de}\cite{bf} A {\em vertex Poisson algebra} is a differential
algebra $(A,\partial)$ equipped with a vertex Lie algebra
structure $(Y_{-},\partial)$ (with the same operator $\partial$)
such that for $a,b,c\in A$,
\begin{equation}\label{gd}
Y_-(a,x)(bc)=(Y_-(a,x)b)c+b(Y_-(a,x)c).
\end{equation}
\end{de}
\begin{rem}
In terms of components, (\ref{gd}) is equivalent to
\begin{equation}\label{hd}
a_n(bc)=(a_nb)c+b(a_nc)\ \ \text{for }a,b,c\in A, n\in \N.
\end{equation}
Here $Y_-(a,x)=\sum_{n\in \N}a_nx^{-n-1}$.
\end{rem}
\begin{coro}\label{dera} For $a\in A$, $n\in \N$,

\begin{enumerate}
\item $a_n$ is derivation of $A$, and

\item $Y_-(a,x)1=0$, and $Y_-(1,x)a=0$. \end{enumerate}
\end{coro}
\begin{proof} The second part of 2. follows from the first part of 2., and (\ref{yhs}).
\end{proof}

A vertex algebra is called {\em commutative} if
$[Y(u,x_1),Y(v,x_2)]=0$ for all $u,v\in V$. (See appendix for the
definition of a vertex algebra.) It was shown in \cite{fhl} that a
vertex algebra $V$ is commutative if and only if $u_nv=0$ for all
$u,v\in V$, $n\in\N$.

\begin{prop}{\cite{bo, bf, li1}} If $V$ is a commutative vertex algebra, then $V$ is a
commutative associative algebra with the product defined by
\begin{equation}
u\cdot v=u_{-1}v \ \ \text{ for } \ \ u,v\in V
\end{equation}
and with ${\bf 1}$ as the identity element. Furthermore, the
operator $\D$ of $V$ is a derivation and $Y(u,x)v=(e^{x\D}u)v\ \
\text{ for }u,v\in V.$

Conversely, for any differential algebra $(A,\partial)$, $(A,Y,1)$
is a commutative vertex algebra where $Y$ is defined by
$Y(a,x)b=(e^{x\partial}a)b$ for $a,b\in A$. This give rise to a
canonical isomorphism between the category of commutative vertex
algebras and the category of differential algebras. Moreover, a
vertex Poisson algebra structure on a vector space $A$ consists of
a commutative vertex algebra structure and a vertex Lie algebra
structure with a compatibility condition.
\end{prop}
\begin{de}\label{vpagrad}\cite{li2} A {\em $\Z$-graded vertex Poisson algebra} is a vertex
Poisson algebra $A$ equipped with a $\Z$-grading $A=\coprod_{n\in
\Z}A_{(n)}$ such that $A$ as an algebra is $\Z$-graded and such
that for $a\in A_{(n)}$, $n,r\in \Z$, $m\in \N$
\begin{eqnarray*}
\partial A_{(r)}&\subset& A_{(r+1)},\\
a_mA_{(r)}&\subset&A_{(r+n-m-1)}.
\end{eqnarray*}
An $\N$-graded vertex Poisson algebra is defined in the obvious
way.
\end{de}
\begin{rem}\label{zcvp} Let $A$ be an $\Z$-graded (respectively, $\N$-graded)
vertex Poisson algebra. Then
\begin{enumerate}
\item $A_{(0)}$ is a commutative associative algebra with the
identity 1. \item For $n\in\Z$ (respectively, $n\in\N$), $A_{(n)}$
is an $A_{(0)}$-module. \item $\partial: A_{(0)}\rightarrow
A_{(1)}$ is a derivation when we consider $A_{(0)}$ as an algebra
and $A_{(1)}$ as an $A_{(0)}$-module.
\end{enumerate}
\end{rem}
\begin{de}\cite{li2} An {\em ideal} of a vertex Poisson algebra $A$ is an ideal
$I$ of $A$ as an associative algebra such that
\begin{eqnarray*}
\partial I&\subset& I\\
a_nI&\subset&I \ \ \text{ for } a\in A, n\in\N.
\end{eqnarray*}
\end{de}
\begin{coro}\cite{li2}\ \
\begin{enumerate}
\item By the half skew symmetry we have $u_nA\subset I$ for $u\in
I$, $n\in\N$.

\item Moreover, the quotient space $A/I$ has a natural vertex
Poisson algebra structure.
\end{enumerate}
\end{coro}

\begin{prop} Let $A$ be a vertex Poisson algebra, and let $I$ be an ideal of $A$ as an associative commutative
algebra. Assume that $\partial I\subset I$ and $I$ is generated by
$W$, that is $I=AW$. If $a_nw\in I$ for all $a\in A, w\in
W,n\in\N$, then I is an ideal of a vertex Poisson algebra $A$.
\end{prop}
\begin{proof} Let $a,a'\in A, w\in W$, $n\in\N$. By (\ref{hd}), we have
$$a_n(a'w)=(a_n a')w+a'(a_nw)\in I.$$ It follows that $I$ is an ideal.
\end{proof}

\begin{prop}\label{rvertexlie}\cite{li2, bf} Let $R$ be a vector space equipped with a linear
operator $\D$ and let $Y^0_-$ be a linear map from $R$ to
$\Hom(R,x^{-1}R[x^{-1}])$. Denote by $S(R)$ the symmetric algebra
over $R$ and we extend $\D$ uniquely to a derivation of $S(R)$.
Then $Y^0_-$ extends to a vertex Poisson algebra structure $Y_-$
on $(S(R),\D)$ if and only if $(R,\D,Y^0_-)$ carries the structure
of a vertex Lie algebra. Furthermore, such an extension is unique.
\end{prop}

\section{From vertex Poisson algebras to Courant algebroids}

First, we review the notions of 1-truncated conformal algebras and
Courant algebroids. Also, we study the relations between these
algebras. Next, we show that for an $\N$-graded vertex Poisson
algebra $V=\coprod_{n\in\N}V_{(n)}$, the homogeneous subspace
$V_{(1)}$ is a Courant $V_{(0)}$-algebroid.

\begin{de}\cite{gms} A {\em 1-truncated conformal algebra} is a graded
vector space $C=C_0\oplus C_1$, equipped with a linear map
$\partial:C_0\rightarrow C_1$ and bilinear operations
$(u,v)\mapsto u_iv$ for $i=0,1$ of degree $-i-1$ on $C=C_0\oplus
C_1$ such that the following axioms hold:

(Derivation) for $a\in C_0$, $u\in C_1$,
\begin{equation}\label{1tcd}
(\partial a)_0=0;\ \ (\partial a)_1=-a_0;\ \
\partial(u_0a)=u_0\partial a
\end{equation}

(Commutativity) for $a\in C_0$, $u,v\in C_1$,
\begin{equation}\label{1tcs}
u_0a=-a_0u;\ \ u_0v=-v_0u+\partial(v_1u);\ \ u_1v=v_1u
\end{equation}

(Associativity) for $\alpha,\beta,\gamma\in C_0\oplus C_1$,
\begin{equation}\label{1tca}
\alpha_0\beta_i\gamma=\beta_i\alpha_0\gamma+(\alpha_0\beta)_i\gamma.
\end{equation}
\end{de}

\begin{rem}\label{vltcf}\ \

\begin{enumerate}

\item Let $A=\coprod_{n\in\N}A_{(n)}$ be an $\N$-graded vertex Lie
algebra. By equations (\ref{hp})-(\ref{ha}), we have that
$A_{(0)}\oplus A_{(1)}$ is a 1-truncated conformal algebra. \item
Consequently, for an $\N$-graded vertex Poisson algebra
$A=\coprod_{n\in\N}A_{(n)}$, $A_{(0)}\oplus A_{(1)}$ is a
1-truncated conformal algebra.
\end{enumerate}
\end{rem}

A {\em Leibniz algebra} is a nonassociative algebra $\Gamma$
satisfying the following condition:
$$u\cdot(v\cdot w)=(u\cdot v)\cdot w+v\cdot(u\cdot w)\ \ \text{ for
}u,v,w\in \Gamma.$$ Any Lie algebra is a Leibniz algebra. In
particular, for any vector space $W$ the general linear Lie
algebra $\mathfrak{gl}(W)$ is a Leibniz algebra. A {\em
representation of Leibniz algebra} $\Gamma$ on a vector space $W$
is a Leibniz algebra homomorphism $\rho$ from $\Gamma$ to
$\mathfrak{gl}(W)$.

Let $\mathcal{A}$ be a unital commutative associative algebra
(over $\C$). A {\em Leibniz $\mathcal{A}$-algebra} is a Leibniz
algebra $\Gamma$ equipped with an $\mathcal{A}$-module structure.
A {\em module} for a Leibniz $\mathcal{A}$-algebra $\Gamma$ is a
vector space $W$ equipped with $\Gamma$-module structure and an
$\mathcal{A}$-module structure.

\begin{prop}\label{1tcf}\cite{gms, liy} Let $C=C_0\oplus C_1$ be a graded vector
space (over $\C$) equipped with a linear map $\partial$ from $C_0$
to $C_1$ and equipped with bilinear maps $(u,v)\mapsto u_iv$ of
degree $-i-1$ on $C=C_0\oplus C_1$ for $i=0,1$. Then $C$ is a
1-truncated conformal algebra if and only if
\begin{enumerate}
\item The pair $(C_1,[\cdot,\cdot])$ carries the structure of a
Leibniz algebra where $[u,v]=u_0v$ for $u,v\in C_1$. \item The
space $C_0$ is a $C_1$-module with $u\cdot a=u_0a$ for $u\in
C_1,a\in C_0$. \item The map $\partial$ is a $C_1$-module
homomorphism. \item The subspace $\partial C_0$ of $C_1$
annihilates the $C_1$-module $C_0\oplus C_1$. \item The bilinear
map $\<\cdot,\cdot\>$ from $C_1\otimes C_1$ to $C_0$ defined by
$\<u,v\>=u_1v$ for $u,v\in C_1$ is a $C_1$-module homomorphism and
furthermore
\begin{eqnarray}
u_0a&=&-a_0u,\\
\<\partial a, u\>&=&-a_0u,\\
{[u,v]+[v,u]}&=&\partial\<u,v\>,\\
\<u,v\>&=&\<v,u\>
\end{eqnarray}
for $a\in C_0,u,v\in C_1$.
\end{enumerate}
\end{prop}

Next, we state the definition of a Courant algebroid. We also
study the relationship between Courant algebroids and 1-truncated
conformal algebras.
\begin{de}\cite{br} Let $\mathcal{A}$ be a unital commutative associative algebra over $\C$. A {\em
Courant $\mathcal{A}$-algebroid} is an $\mathcal{A}$-module
$\mathcal{B}$ equipped with
\begin{enumerate}
\item a structure of Leibniz algebra $[\ \ , \ \
]:\mathcal{B}\otimes_{\C} \mathcal{B}\rightarrow \mathcal{B}$,

\item a homomorphism of Leibniz $\mathcal{A}$-algebras $\pi:
\mathcal{B} \rightarrow Der(\mathcal{A})$,

\item a symmetric $\mathcal{A}$-bilinear pairing $\<\ \ , \ \ \>:
\mathcal{B}\otimes_{\mathcal{A}} \mathcal{B}\rightarrow
\mathcal{A}$,

\item a derivation $\partial :\mathcal{A}\rightarrow \mathcal{B}$
such that $\pi\circ
\partial =0$ which satisfy
\begin{eqnarray}
& &[u,av]=a[u,v]+\pi(u)(a)v,\label{c1}\\
& &\<[u,v],w\>+\<v,[u,w]\>=\pi(u)\<v,w\>,\label{c2}\\
& &[u,\partial a]=\partial(\pi(u)a),\label{c3}\\
& &\<u,\partial a\>=\pi(u)a,\label{c4}\\
& &[u,v]+[v,u]=\partial(\<u,v\>)\label{c5},
\end{eqnarray}
for $a\in \mathcal{A}, u,v,w\in \mathcal{B}$.
\end{enumerate}
\end{de}
\begin{coro}\label{hom}\ \

\begin{enumerate}
\item $\partial (\mathcal{A})$ annihilates $\mathcal{A}$ and
$\mathcal{B}$. \item $\< \ \ ,\ \ \>$ and $\partial$ are
$\mathcal{B}$-module homomorphisms.
\end{enumerate}
\end{coro}

\begin{proof}

For 1., since $\pi\circ\partial=0$, it follows immediately that
$\partial (\mathcal{A})$ annihilates $\mathcal{A}$. Next, we will
show that $\partial(\mathcal{A})$ annihilates $\mathcal{B}$. Let
$a\in \mathcal{A}$, $u\in \mathcal{B}$. By (\ref{c3})-(\ref{c5}),
we have $$[\partial (a),u]=-[u,\partial (a)]+\partial(\<\partial
(a),u\>)=-\partial (\pi(u)a)+\partial(\pi(u)a)=0.$$

2. follows immediately from (\ref{c2}), (\ref{c3}).
\end{proof}

By the definition of a Courant algebroid, Proposition \ref{1tcf},
and Corollary \ref{hom}, we have the following.
\begin{prop}\label{cor1trun} Let $\mathcal{A}$ be a unital commutative associative
algebra, and let $\mathcal{B}$ be an $\mathcal{A}$-module. Let
$\partial:\mathcal{A}\rightarrow \mathcal{B}$ be a derivation.
Then a Courant $\mathcal{A}$-algebroid structure on $\mathcal{B}$
is exactly equivalent to a 1-truncated conformal algebra structure
on $C=\mathcal{A}\oplus \mathcal{B}$ with
\begin{eqnarray*}
& & a_ia'=0,\\
& &u_0v=[u,v],\ \ u_1v=\<u,v\>,\\
& &u_0a=\pi(u)(a),\ \ a_0u=-u_0a
\end{eqnarray*}
for $a,a'\in \mathcal{A}$, $u,v\in \mathcal{B}$, $i=0,1$ such that
\begin{eqnarray}
& &(au)_0a'=a(u_0a'),\label{dera1}\\
& &(au)_1v=a(u_1v)=u_1(av),\label{syma}\\
& &u_0(av)=a(u_0v)+(u_0a)v,\label{dera2}\\
& &u_0(aa')=a(u_0a')+(u_0a)a'.\label{dec}
\end{eqnarray}
\end{prop}

\begin{coro}\label{dee} Let $e$ be the identity of $\mathcal{A}$. Then for $u\in \mathcal{B}$, $u_0e=0$.
\end{coro}
\begin{proof} This follows immediately from (\ref{dec}).
\end{proof}

Let $(A=\coprod_{n\in\Z}A_{(n)},\partial)$ be an $\N$-graded
vertex Poisson algebra. By Remark \ref{zcvp}, we have $A_{(0)}$ is
a commutative associative algebra with the identity 1 and
$A_{(1)}$ is an $A_{(0)}$-module. Moreover,
$\partial:A_{(0)}\rightarrow A_{(1)}$ is a derivation. By Remark
\ref{vltcf}, we also have that $A_{(0)}\oplus A_{(1)}$ is a
1-truncated conformal algebra.

\begin{thm} $A_{(1)}$ is, in fact, a Courant $A_{(0)}$-algebroid.
\end{thm}
\begin{proof} By Proposition \ref{cor1trun}, it is enough to show that
(\ref{dera1})-(\ref{dec}) hold on $A_{(0)}\oplus A_{(1)}$. Let
$a,a'\in A_{(0)}, u,v\in A_{(1)}$. By (\ref{hd}), we have
\begin{eqnarray*}
& &u_0(av)=a(u_0v)+(u_0a)v,\\
& &u_0(aa')=(u_0a)a'+a(u_0a').
\end{eqnarray*}
Hence, (\ref{dera2}), (\ref{dec}) hold.

Let $a\in A_{(0)}, u,v\in A_{(1)}$. By (\ref{hd}), (\ref{1tcs}),
we have
\begin{eqnarray*}
& &u_1(av)=a(u_1v), \text{ and }\\
& &(au)_1v=v_1(au)=a(v_1u)=a(u_1v). \end{eqnarray*} Therefore,
(\ref{syma}) holds.

Next, we will show that (\ref{dera1}) holds. Let $a,a'\in
A_{(0)},u\in A_{(1)}$. By (\ref{hs}), (\ref{hd}), (\ref{1tcs}), we
have
\begin{eqnarray*}
(au)_0a'
&=&\sum_{i\geq 0}(-1)^{i+1}\frac{1}{i!}\partial ^i(a'_i(au))\\
&=&-a'_0(au)\\
&=&-\{(a'_0a)u+a(a'_0u)\}\\
&=&-a(-u_0a')\\
&=&a(u_0a').
\end{eqnarray*}
Therefore, $A_{(1)}$ is a Courant $A_{(0)}$-algebroid.
\end{proof}

\section{Vertex Poisson algebras associated with 1-truncated conformal algebras}

In this section, we construct an $\N$-graded vertex Poisson
algebra from any 1-truncated conformal algebra.

Let $(\mathbb{A}\oplus \mathbb{B},\partial)$ be a 1-truncated
conformal algebra. We set
$$C(\mathbb{A}\oplus \mathbb{B})=\C[\D]\otimes_{\C}(\mathbb{A}\oplus\mathbb{B}),$$ where $\D$ is a formal variable.
We define subspaces $C(\mathbb{A})$ and $C(\mathbb{B})$ in the
obvious way. Additionally, we consider
$\mathbb{A}\oplus\mathbb{B}$ as a subspace of
$C(\mathbb{A}\oplus\mathbb{B})$ via the following map
$$\mathbb{A}\oplus\mathbb{B}\rightarrow C(\mathbb{A}\oplus\mathbb{B}),\ \ a+b\mapsto 1\otimes a+1\otimes b.$$
Also, we define a linear operator $\D$ on
$C(\mathbb{A}\oplus\mathbb{B})$ by
\begin{equation}\label{eqd}\D(\D^n\otimes u)=\D^{n+1}\otimes u\text{ for }
u\in\mathbb{A}\oplus \mathbb{B}, n\in\N,\end{equation} and set
\begin{equation}\label{hatpartial}
\hat{\partial}=1\otimes \partial-\D \otimes 1:
{C}(\mathbb{A})\rightarrow {C}(\mathbb{A}\oplus \mathbb{B}).
\end{equation}

Next, we define
\begin{eqnarray*}
& & \deg(\D ^n\otimes a)=n\ \ \text{ for }a\in \mathbb{A},n\in\N,\\
& & \deg(\D ^n\otimes b)=n +1\ \ \text{ for } b\in \mathbb{B},
n\in\N.
\end{eqnarray*}
Then $C(\mathbb{A}\oplus \mathbb{B})$ becomes an $\N$-graded
vector space:
$$C(\mathbb{A}\oplus \mathbb{B})=\coprod_{n\in\N}C(\mathbb{A}\oplus \mathbb{B})_{(n)},$$
where $C(\mathbb{A})_{(0)}=\mathbb{A}$, and for $n\geq 1$
\begin{eqnarray}
C(\mathbb{A}\oplus\mathbb{B})_{(n)}&=&\D ^n\otimes \mathbb{A}\oplus \D^{n-1}\otimes \mathbb{B}\nonumber\\
&=&\{\ \ \D ^n\otimes a,\ \ \D ^{n-1}\otimes b\ \ |\ \
a\in\mathbb{A},\ \ b\in\mathbb{B}\ \ \}\label{gradab}.
\end{eqnarray}
The linear map $\hat{\partial}$ is a homogeneous of degree 1 and
for $n\geq 1$, we have
$$(\hat{\partial}C(\mathbb{A}))_{(n)}=\{\ \ \hat{\partial}(\D ^{n-1}\otimes a)\ \ |\ \ a\in \mathbb{A}\ \ \}.$$ Moreover,
$\D(\hat{\partial}C(\mathbb{A}))\subset\hat{\partial}C(\mathbb{A})$.

We define a linear map
\begin{eqnarray*}
Y^0_-:C(\mathbb{A}\oplus \mathbb{B})&\rightarrow&
\Hom(C(\mathbb{A}\oplus \mathbb{B}),x^{-1}C(\mathbb{A}\oplus\mathbb{B})[x^{-1}])\\
u&\mapsto&Y^0_-(u,x)=\sum_{n\in\N}u_nx^{-n-1} \end{eqnarray*} in
the following way: for $a,a'\in\mathbb{A}$, $b,b'\in\mathbb{B}$,
$u,v\in\mathbb{A}\oplus\mathbb{B}, n\geq 1$, $m\in\N$,
\begin{eqnarray}
& &Y^0_{-}(a,x)a'=0,\label{yaa'}\\
& &Y^0_{-}(a,x)b=a_0bx^{-1},\label{yab}\\
& &Y^0_{-}(b,x)a=b_0ax^{-1},\label{yba}\\
& &Y^0_{-}(b,x)b'=b_0b'x^{-1}+b_1b'x^{-2},\label{ybb'}\\
&
&Y^0_-(u,x)e^{x_1\D}v=e^{x_1\D}e^{-x_1\frac{d}{dx}}Y^0_-(u,x)v,\label{uedv}\\
& &Y^0_-(\D ^n \otimes u,x)\D ^m \otimes
v=Sing(e^{x\D}\left(-\frac{d}{dx}\right)^mY^0_-(v,-x)\D^n\otimes
u).\label{ydu}
\end{eqnarray}
\begin{rem} Notice that (\ref{uedv}) is equivalent to
\begin{equation}\label{yudv}
Y^0_{-}(u,x)\D ^n \otimes
v=\sum_{i=0}^n\frac{n!}{(n-i)!i!}(-1)^{n-i}\D
^i\left(\frac{d}{dx}\right)^{n-i}Y^0_-(u,x)v
\end{equation} for $u,v\in \mathbb{A}\oplus \mathbb{B}$, $n\in\N$.
\end{rem}

The following proposition will play an important role for the rest
of this section.
\begin{prop}\label{grad1} For $a,a'\in\mathbb{A}$, $b,b'\in\mathbb{B}$, $n\in\N$,
we have
\begin{eqnarray}
Y^0_-(a,x)\D ^n\otimes a'&=&0\label{yada'}\\
Y^0_-(a,x)\D ^n \otimes b&=&\sum_{i=0}^n\frac{n!}{(n-i)!}\D
^{n-i}\otimes (a_0b)x^{-1-i}\label{yadb}\\
Y^0_-(b,x)\D ^n \otimes a&=&\sum_{i=0}^n\frac{n!}{(n-i)!}\D
^{n-i}\otimes (b_0a)x^{-1-i}\label{ybda}\\
Y^0_-(b,x)\D ^n \otimes b'&=&\sum_{i=0}^n\frac{n!}{(n-i)!}\D
^{n-i}\otimes (b_0b')x^{-1-i}\nonumber\\
& &+\sum_{i=1}^{n+1} \frac{in!}{(n-i+1)!}\D^{n-i+1}\otimes
(b_1b')x^{-i-1}\label{ybdb'}
\end{eqnarray}
\end{prop}
\begin{proof} It follows immediately from (\ref{yaa'})-(\ref{ybb'}),
(\ref{yudv}).\end{proof}

Next, we will show that $(C(\mathbb{A}\oplus
\mathbb{B})/\hat{\partial}C(\mathbb{A}), Y^0_-,\D)$ is a vertex
Lie algebra.
\begin{prop}\label{acideal} The subspace $\hat{\partial}C(\mathbb{A})$ of the
nonassociative algebra $C(\mathbb{A}\oplus\mathbb{B})$ is a
2-sided ideal.
\end{prop}
\begin{proof} First we show that for $u\in
\mathbb{A}\oplus \mathbb{B}$, $a\in\mathbb{A}$, $ m,n\in\N$,
$$Y^0_-(\D^n \otimes u,x)\hat{\partial}(\D ^m\otimes a)\in
x^{-1}\hat{\partial}C(\mathbb{A})[x^{-1}].$$ Let
$a,a'\in\mathbb{A}$, $m\in\N$, $n\geq 1$. By
(\ref{1tcd})-(\ref{1tcs}), (\ref{ydu}),
(\ref{yada'})-(\ref{ybda}), we have
$$Y^0_-(a',x)\hat{\partial}(\D^m\otimes a)=Y^0_-(a',x)\left(\D^m\otimes \partial(a)-\D
^{m+1}\otimes a\right)=0,$$ and
\begin{eqnarray*}
& &Y^0_-(\D ^n\otimes a',x)\hat{\partial}(\D^m\otimes a)\\
&=&Y^0_-(\D^n\otimes a',x)\left(\D^m\otimes \partial(a)-\D
^{m+1}\otimes a\right)\\
&=&Sing(e^{x\D}\left(-\frac{d}{dx}\right)^{m}Y^0_-(\partial(a),x)D^n
\otimes a')\\
& &-Sing(e^{x\D}\left(-\frac{d}{dx}\right)^{m+1}Y^0_-(a,-x)\D^n\otimes a')\\
&=&0.
\end{eqnarray*}
Let $a\in\mathbb{A}$, $b\in\mathbb{B}$, $m\in\N$, $n\geq 1$. By
(\ref{1tcd})-(\ref{1tcs}), (\ref{ydu}),
(\ref{yadb})-(\ref{ybdb'}), we have
\begin{eqnarray*}
& &Y^0_-(b,x)\hat{\partial}(\D^m\otimes a)\\
&=&Y^0_-(b,x)\D ^m\otimes \partial (a)-Y^0_-(b,x)\D ^{m+1}\otimes a\\
&=&\sum_{i=0}^m\frac{m!}{(m-i)!}\D ^{m-i}\otimes
b_0\partial(a)x^{-1-i}+\sum_{i=1}^{m+1}\frac{im!}{(m-i+1)!}\D
^{m-i+1}\otimes b_1\partial(a)x^{-i-1}\\
& &-\sum_{i=0}^{m+1}\frac{(m+1)!}{(m+1-i)!}\D ^{m+1-i}\otimes
b_0a x^{-1-i}\\
&=&\sum_{i=0}^m\frac{m!}{(m-i)!}\D ^{m-i}\otimes
\partial((\partial a)_1b)x^{-1-i}+\sum_{i=1}^{m+1}\frac{im!}{(m-i+1)!}\D
^{m-i+1}\otimes (\partial a)_1bx^{-i-1}\\
& &-\sum_{i=0}^{m+1}\frac{(m+1)!}{(m+1-i)!}\D ^{m+1-i}\otimes
(\partial a)_1b x^{-1-i}\\
&=&\sum_{i=0}^m\frac{m!}{(m-i)!}\D ^{m-i}\otimes
\partial((\partial a)_1b)x^{-1-i}-\sum_{i=0}^{m}\frac{m!}{(m-i)!}\D
^{m+1-i}\otimes
(\partial a)_1bx^{-1-i}\\
&=&\sum_{i=0}^m\frac{m!}{(m-i)!}\hat{\partial}\left(\D
^{m-i}\otimes (\partial a)_1b \right)x^{-1-i},
\end{eqnarray*}
and
\begin{eqnarray*}
& &Y^0_-(\D ^n\otimes b,x)\hat{\partial}(\D^m\otimes a)\\
&=&Y^0_-(\D ^n\otimes b,x)\D ^m\otimes \partial (a)-Y^0_-(\D
^n\otimes b,x)\D ^{m+1}\otimes a\\
&=&Sing(e^{x\D}\left(-\frac{d}{dx}\right)^mY^0_-(\partial
(a),-x)\D ^n\otimes b)\\
& &- Sing (e^{x\D}\left(-\frac{d}{dx}\right)^{m+1}Y^0_-(a,-x)\D ^n
\otimes b)\\
&=&Sing (e^{x\D}\left(-\frac{d}{dx}\right)^m\sum_{i=0}^n
\frac{n!}{(n-i)!}\D ^{n-i}\otimes (\partial(a)_0b)(-x)^{-i-1})\\
& &+Sing(
e^{x\D}\left(-\frac{d}{dx}\right)^{m}\sum_{i=1}^{n+1}\frac{in!}{(n-i+1)!}\D^{n-i+1}\otimes
\partial (a)_1b(-x)^{-i-1})\\
& &-Sing(
e^{x\D}\left(-\frac{d}{dx}\right)^{m+1}\sum_{i=0}^{n}\frac{n!}{(n-i)!}\D^{n-i}\otimes
(a_0b)(-x)^{-1-i})\\
&=&0.
\end{eqnarray*}
These imply that for $u\in \mathbb{A}\oplus \mathbb{B}$,
$a\in\mathbb{A}$, $m,n\in\N$, $$Y^0_-(\D^n \otimes
u,x)\hat{\partial}(\D ^m\otimes a)\in
x^{-1}\hat{\partial}C(\mathbb{A})[x^{-1}].$$

Next, we will show that for $u\in \mathbb{A}\oplus\mathbb{B}$,
$a\in\mathbb{A}$, $m,n\in\N$,
$$Y^0_-(\hat{\partial}(\D ^m\otimes a),x)\D ^n\otimes u\in
x^{-1}\hat{\partial}C(\mathbb{A})[x^{-1}].$$ Let
$a,a'\in\mathbb{A}$, $b\in\mathbb{B}$, $n\in\N$. By (\ref{1tcd}),
(\ref{ydu}), (\ref{yada'}), (\ref{ybda})-(\ref{ybdb'}), we have
\begin{eqnarray*}
Y^0_-(\hat{\partial}(a),x)\D^n\otimes a'
&=&Y^0_-(\partial(a),x)\D^n\otimes a'-Y(\D\otimes a,x)\D^n\otimes
a'\\
&=&-Sing(e^{x\D}\left(-\frac{d}{dx}\right)^nY^0_-(a',-x)\D \otimes
a)\\
&=&0,
\end{eqnarray*}
and
\begin{eqnarray*}
& &Y^0_-(\hat{\partial}(1\otimes a),x)\D^n\otimes b\\
&=&Y^0_-(1\otimes
\partial a,x)\D^n\otimes b-Y^0_-(\D\otimes
a,x)\D^n\otimes b\\
&=&\sum_{i=0}^n\frac{n!}{(n-i)!}\D^{n-i}\otimes (\partial
a)_{(0)}bx^{-1-i}+\sum_{i=1}^{n+1}\frac{in!}{(n-i+1)!}\D^{n-i+1}\otimes
(\partial a)_1bx^{-i-1}\\
& &-Sing(e^{x\D}\left(-\frac{d}{dx}\right)^nY^0_-(b,-x)\D\otimes
a)\\
&=&-\sum_{i=1}^{n+1}\frac{in!}{(n-i+1)!}\D^{n-i+1}\otimes
a_0bx^{-i-1}\\
&
&+Sing(e^{x\D}\left(-\frac{d}{dx}\right)^n(\D\otimes a_0b(-x)^{-1}+1\otimes a_0b(-x)^{-2}))\\
&=&-\sum_{i=1}^{n+1}\frac{in!}{(n-i+1)!}\D^{n-i+1}\otimes
a_0bx^{-i-1}-\sum_{j=0}^n\frac{n!}{j!}\D^{j+1}\otimes
a_0bx^{-1-n+j}\\
& &+\sum_{j=0}^{n+1}\frac{(n+1)!}{j!}\D^j\otimes a_0bx^{-2-n+j}\\
&=&-\sum_{i=1}^{n+1}\frac{in!}{(n-i+1)!}\D^{n-i+1}\otimes
a_0bx^{-i-1}-\sum_{i=0}^n\frac{n!}{(n-i)!}\D^{n-i+1}\otimes
a_0bx^{-1-i}\\
& &+\sum_{i=0}^{n+1}\frac{(n+1)!}{(n-i+1)!}\D^{n-i+1}\otimes a_0bx^{-1-i}\\
&=&0.
\end{eqnarray*}
Let $a\in\mathbb{A}$, $b\in\mathbb{B}$, $n\in\N$, $m\geq 1$. By
(\ref{1tcd})-(\ref{1tcs}), (\ref{ydu}),
(\ref{yada'})-(\ref{yadb}), we have
\begin{eqnarray*}
Y^0_-(\hat{\partial}(\D^m\otimes a),x)\D ^n\otimes a'
&=&Y^0_-(\D^m\otimes \partial a,x)\D ^n\otimes a'-Y(\D
^{m+1}\otimes a,x)\D^n\otimes a'\\
&=&Sing
(e^{x\D}\left(-\frac{d}{dx}\right)^nY^0_-(a',-x)\D^m\otimes
\partial a)\\
&=&0.
\end{eqnarray*}
Since $\D(\hat{\partial}(C(\mathbb{A}))\subset
\hat{\partial}(C(\mathbb{A}))$ and
$Y^0_-(b,-x)\hat{\partial}(C(\mathbb{A}))\subset
x^{-1}\hat{\partial}(C(\mathbb{A}))[x^{-1}]$, we have
\begin{eqnarray*}
& &Y^0_-(\hat{\partial}(\D^m\otimes a),x)\D ^n\otimes b\\
&=&Y^0_-(\D ^m\otimes \partial (a),x)\D ^n\otimes b-Y^0_-(\D
^{m+1}\otimes a,x)\D^n\otimes b\\
&=&Sing(e^{x\D}\left(-\frac{d}{dx}\right)^nY^0_-(b,-x)\D^m\otimes
\partial
(a))\\
&
&-Sing(e^{x\D}\left(-\frac{d}{dx}\right)^nY^0_-(b,-x)\D^{m+1}\otimes
a)\\
&=&Sing(e^{x\D}\left(-\frac{d}{dx}\right)^nY^0_-(b,-x)\hat{\partial}(\D
^m\otimes a))\in x^{-1}\hat{\partial}(C(\mathbb{A}))[x^{-1}].
\end{eqnarray*}
It follows that for $u\in \mathbb{A}\oplus\mathbb{B}$,
$a\in\mathbb{A}$, $m,n\in\N$, we have
$$Y^0_-(\hat{\partial}(\D ^m\otimes a),x)\D ^n\otimes u\in
x^{-1}\hat{\partial}C(\mathbb{A})[x^{-1}].$$ Hence, $\hat\partial
C(\mathbb{A})$ is a 2-sided ideal of
$C(\mathbb{A}\oplus\mathbb{B})$.
\end{proof}

Set
\begin{equation}\label{vlc}
\mathcal{C}=C(\mathbb{A}\oplus
\mathbb{B})/\hat{\partial}C(\mathbb{A}).\end{equation} Let $\rho$
be a projection map:
$$C(\mathbb{A}\oplus\mathbb{B})\rightarrow \mathcal{C};\D
^n\otimes u\mapsto \D ^n \otimes u +\hat{\partial}C(\mathbb{A}).$$
For $u\in \mathbb{A}\oplus \mathbb{B}, n\in\N$, we set $\D^n
u=\rho (\D ^n\otimes u).$

By (\ref{hatpartial}), (\ref{gradab}) and (\ref{vlc}), we have the
following.
\begin{prop}\label{Dpartial} For $a\in \mathbb{A}$, $n\in\N$,
$\D ^n(\partial a)=\D ^{n+1}a.$ Moreover, for $n\geq 1$,
$\mathcal{C}_{(n)}=\{\D^{n-1}(b)\ \ |\ \ b\in \mathbb{B}\}.$
\end{prop}
\begin{prop}\label{grading1} For $u\in\mathbb{A}\oplus \mathbb{B}$, $i\geq 0$
$n\in\N$, $u_i\mathcal{C}_{(n)}\subset\mathcal{C}_{(n+\deg u
-i-1)}$. \end{prop}
\begin{proof} This follows immediately from Propositions
\ref{grad1}, \ref{acideal}.
\end{proof}

Next, we will show that $(\mathcal{C}, Y^0_-,\D)$ is a vertex Lie
algebra.
\begin{lem} \label{sing} For $u,v\in\mathbb{A}\oplus \mathbb{B}$, $n\in\N$, we have
$$Y^0_-(u,x)\D ^n
v=Sing(e^{x\D}\left(-\frac{d}{dx}\right)^nY^0_-(v,-x)u).$$
\end{lem}
\begin{proof} We first show that for $u,v\in\mathbb{A}\oplus \mathbb{B}$,
$$Y^0_-(u,x)v=Sing(e^{x\D}Y^0_-(v,-x)u).$$
Let $a,a'\in\mathbb{A}$, $b,b'\in\mathbb{B}$. By (\ref{1tcs}),
(\ref{yaa'})-(\ref{ybb'}), and Proposition \ref{Dpartial}, we have
\begin{eqnarray*}
& &Y^0_{-}(a,x)a'=0=Sing(e^{x\D}Y^0_{-}(a',-x)a),\\
& &Y^0_{-}(a,x)b=a_0bx^{-1}=-b_0ax^{-1}=Sing
(e^{x\D}Y^0_{-}(b,-x)a),\\
&
&Y^0_{-}(b,x)a=b_0ax^{-1}=-a_0bx^{-1}=Sing(e^{x\D}Y^0_{-}(a,-x)b),
\end{eqnarray*}
and
\begin{eqnarray*}
& &Sing(e^{x\D}Y^0_{-}(b',-x)b)\\
&=&Sing(e^{x\D}(-b'_0bx^{-1}+b'_1bx^{-2}))\\
&=&-b'_0bx^{-1}+b'_1bx^{-2}+\D(b'_1b)x^{-1}\\
&=&-b'_0bx^{-1}+b'_1bx^{-2}+\partial(b'_1b)x^{-1}\\
&=&b_0b'x^{-1}+b_1b'x^{-2}\\
&=&Y^0_{-}(b,x)b'.\end{eqnarray*} These imply that for
$u,v\in\mathbb{A}\oplus \mathbb{B}$,
\begin{equation}\label{hsuv}
Y^0_{-}(u,x)v=Sing (e^{x\D }Y^0_{-}(v,-x)u).
\end{equation}

Next, we show that for $u,v\in\mathbb{A}\oplus\mathbb{B}$, $n\geq
1$, $$Y^0_-(u,x)\D ^n v=Sing(e^{x\D}\left(-\frac{d}{dx}\right)^n
Y^0_-(v,-x)u).$$ Clearly, for all $a,a'\in \mathbb{A}$, $n\geq 1$,
we have $$Y^0_-(a,x)\D
^na'=0=Sing(e^{x\D}\left(-\frac{d}{dx}\right)^nY^0_-(a',x)a).$$
Let $a\in \mathbb{A}$, $b\in\mathbb{B}$, $n\geq 1$. By
(\ref{1tcs}), (\ref{yba}), (\ref{yadb}), we have
\begin{eqnarray*}
& &Sing(e^{x\D}\left(-\frac{d}{dx}\right)^nY^0_-(b,-x)a)\\
&=&Sing(e^{x\D}\left(-\frac{d}{dx}\right)^n b_0a(-x)^{-1})\\
&=&-\sum_{i=0}^n\frac{n!}{i!}\D ^i b_0ax^{-1-n+i}\\
&=&-\sum_{i=0}^n\frac{n!}{(n-i)!}\D^{n-i}b_0ax^{-1-i}\\
&=&Y^0(a,x)\D^n b.
\end{eqnarray*}
Similarly, we have $$Y^0_-(b,x)\D
^na=Sing(e^{x\D}\left(-\frac{d}{dx}\right)^nY^0_-(a,-x)b).$$ Let
$b,b'\in\mathbb{B}$, $n\geq 1$. By (\ref{1tcs}), (\ref{ybb'}),
(\ref{ybdb'}), and Proposition \ref{Dpartial}, we have
\begin{eqnarray*}
& &Sing(e^{x\D}\left(-\frac{d}{dx}\right)^nY^0_-(b,-x)b')\\
&=&Sing(e^{x\D}\left(-\frac{d}{dx}\right)^n\{b_0b'(-x)^{-1}+b_1b'(-x)^{-2}\})\\
&=&\sum_{i=0}^n\frac{n!}{i!}\D
^i(-b_0b')x^{-1-n+i}+\sum_{i=0}^{n+1}\frac{(n+1)!}{i!}\D^ib_1b'x^{-2-n+i}\\
&=&\sum_{i=0}^n\frac{n!}{i!}\D
^i(b'_0b-\partial(b'_1b))x^{-1-n+i}+\sum_{i=0}^{n+1}\frac{(n+1)!}{i!}\D^ib_1b'x^{-2-n+i}\\
&=&\sum_{i=0}^n\frac{n!}{i!}(\D^ib'_0b-\D^{i+1}b'_1b)x^{-1-n+i}+\sum_{i=0}^{n+1}\frac{(n+1)!}{i!}\D^ib_1b'x^{-2-n+i}\\
&=&\sum_{j=0}^n\frac{n!}{(n-j)!}(\D^{n-j}b'_0b-\D^{n-j+1}b'_1b)x^{-1-j}+\sum_{j=0}^{n+1}\frac{(n+1)!}{(n+1-j)!}\D^{n+1-j}b_1b'x^{-1-j}\\
&=&Y^0_-(b',x)\D^n b.
\end{eqnarray*}
Hence, we can conclude that $Y^0_-(u,x)\D ^n
v=Sing(e^{x\D}\left(-\frac{d}{dx}\right)^nY^0_-(v,-x)u)$ for all
$u,v\in\mathbb{A}\oplus\mathbb{B}$, $n\in\N$.
\end{proof}
\begin{prop}\label{bracdy} For $u,v\in\mathbb{A}\oplus \mathbb{B}, n\in\N$,
\begin{equation}\label{bracued}
[\D, Y^0_-(\D ^n
u,x)]e^{x_1\D}v=\frac{d}{dx}Y^0_-(\D^nu,x)e^{x_1\D}v=Y^0_-(\D(\D^n
u),x)e^{x_1\D}v.
\end{equation}
\end{prop}
\begin{proof} We will follow the proof of Proposition 3.10 in \cite{li2} very
closely. First we show that for $u,v\in\mathbb{A}\oplus
\mathbb{B}$, $n\in\N$,
$$[\D,
Y^0_-(\D^nu,x)]e^{x_1\D}v=\frac{d}{dx}Y^0_-(\D^nu,x)e^{x_1\D}v.$$
If we differentiate (\ref{uedv}) with respect to $x_1$, then we
have
\begin{eqnarray*}
Y^0_-(u,x)\D e^{x_1\D}v&=&\D e^{x_1\D}e^{-x_1\frac{d}{dx}}Y^0_-(u,x)v-e^{x_1\D}\frac{d}{dx}e^{-x_1\frac{d}{dx}}Y^0_-(u,x)v\\
&=&\D Y^0_-(u,x)e^{x_1\D}v-\frac{d}{dx}Y^0_-(u,x)e^{x_1\D}v \ \
(\text{by } (\ref{uedv})).
\end{eqnarray*}
Hence, for $u,v\in\mathbb{A}\oplus \mathbb{B}$,
\begin{equation}\label{dydx}[\D,
Y^0_-(u,x)]e^{x_1\D}v=\frac{d}{dx}Y^0_-(u,x)e^{x_1\D}v.\end{equation}
Let $u,v\in\mathbb{A}\oplus\mathbb{B}$, $m\in\N$, $n\geq 1$. By
(\ref{dxsing}), (\ref{ydu}), we have
\begin{eqnarray*}
& &[\D, Y^0_-(\D^n u,x)]\D^mv\\
&=&\D Y^0_-(\D ^n u,x)\D^mv-Y^0_-(\D ^n u,x)\D^{m+1}v\\
&=&Sing(\D e^{x\D}\left(-\frac{d}{dx}\right)^mY^0_-(v,-x)\D^n
u)-Sing(e^{x\D}\left(-\frac{d}{dx}\right)^{m+1}Y^0_-(v,-x)\D^n
u)\\
&=&\frac{d}{dx}Sing(e^{x\D}\left(-\frac{d}{dx}\right)^mY^0_-(v,-x)\D^n
u)\\
&=&\frac{d}{dx}Y^0_-(\D^nu,x)\D^mv. \end{eqnarray*} Hence, for
$u,v\in\mathbb{A}\oplus\mathbb{B}$, $n\in\N$, $[\D, Y^0_-(\D ^n
u,x)]e^{x_1\D}v=\frac{d}{dx}Y^0_-(\D ^n u,x)e^{x_1\D}v$.

Next, we will show that for $u,v\in\mathbb{A}\oplus \mathbb{B}$,
$n\in\N$, $$Y^0_-(\D (\D ^n u),x)e^{x_1\D}v=\frac{d}{dx}Y^0_-(\D
^n u,x)e^{x_1\D}v.$$ By (\ref{ydu}), (\ref{dydx}), and Lemma
\ref{sing}, we have that for $u,v\in\mathbb{A}\oplus\mathbb{B},
n,m\in\N$,
\begin{eqnarray*}
Y^0_-(\D (\D^n u),x)\D ^m
v&=&Sing(e^{x\D}\left(-\frac{d}{dx}\right)^mY^0_-(v,-x)\D ^{n+1}u)\\
&=&Sing(e^{x\D}\left(-\frac{d}{dx}\right)^m(\D
Y^0_-(v,-x)+\frac{d}{dx}Y^0_-(v,-x))\D^{n}u)\\
&=&Sing(e^{x\D}\D\left(-\frac{d}{dx}\right)^mY^0_-(v,-x)\D^nu)\\
& &-Sing(e^{x\D}\left(-\frac{d}{dx}\right)^{m+1}Y^0_-(v,-x)\D^n
u)\\
&=&\frac{d}{dx}Sing(e^{x\D}\left(-\frac{d}{dx}\right)^mY^0_-(v,-x)\D^n
u)\\
&=&\frac{d}{dx}Y^0_-(\D ^n u,x)\D^m v.
\end{eqnarray*}
Therefore, for $u,v\in\mathbb{A}\oplus \mathbb{B}, n\in\N$,
$$[\D, Y^0_-(\D ^n
u,x)]e^{x_1\D}v=\frac{d}{dx}Y^0_-(\D^nu,x)e^{x_1\D}v=Y^0_-(\D
(\D^n u),x)e^{x_1\D}v.$$
\end{proof}

\begin{coro}\label{exd} For any $u\in\mathbb{A}\oplus\mathbb{B}$, $n\in\N$,
$$e^{x_1\D}Y^0_-(\D ^nu,x)e^{-x_1\D}=e^{x_1\frac{d}{dx}}Y^0_-(\D ^nu,x)=Y^0_-(e^{x_1\D}\D ^nu,x)=Y^0_-(\D ^nu,x+x_1).$$
\end{coro}
\begin{coro}\label{dui} For $u\in\mathbb{A}\oplus\mathbb{B}$, $n\geq
1$,
$$(D^n u)_i=\left\{%
\begin{array}{ll}
    0, & \hbox{if \ \ $0\leq i< n$;} \\
    (-1)^n\frac{i!}{(i-n)!}u_{i-n}, & \hbox{if \ \ $i\geq n$.} \\
\end{array}%
\right.$$ Furthermore, we have
$v_i\mathcal{C}_{(n)}\subset\mathcal{C}_{(n+m-i-1)}$ for all $v\in
\mathcal{C}_{(m)}$, $i\in\N$.
\end{coro}
\begin{proof} It follows immediately from Proposition \ref{grading1} and
(\ref{bracued}).\end{proof}
\begin{coro}\label{lhs} For $u,v\in \mathbb{A}\oplus \mathbb{B}$, $n,m\in\N$,
$$Y^0_-(\D ^n u,x)\D^m v=Sing(e^{x\D}Y^0_-(\D^m v,-x)\D ^n u).$$
\end{coro}
\begin{proof} It follows from Lemma \ref{sing}, (\ref{ydu}),
(\ref{bracued}).\end{proof}

\begin{prop}\label{bskew} For $u,v\in\mathbb{A}\oplus\mathbb{B}$, $m,n\in\N$,
$$[Y^0_-(\D^nu,x_1),Y^0_-(\D^m v,x_2)]=Sing(Y^0_-(Y^0_-(\D^n
u,x_1-x_2)\D^m v,x_2)).$$
\end{prop}
\begin{proof} We first show that for $u,v,w\in \mathbb{A}\oplus
\mathbb{B}$,
\begin{eqnarray*}
&&{Y}^0_-(u,x_1)Y^0_-(v,x_2)w-{Y}^0_-(v,x_2)Y^0_-(u,x_1)w\nonumber\\
&=&Sing (Y^0_-(Y^0_-(u,x_1-x_2)v,x_2)w).
\end{eqnarray*}
By Proposition \ref{haf}, it is enough to show that for $u,v,w\in
\mathbb{A}\oplus \mathbb{B}$ from an ordered basis of
$\mathbb{A}\oplus \mathbb{B}$ with $u\leq v\leq w$,
\begin{eqnarray}
&&{Y}^0_-(u,x_1)Y^0_-(v,x_2)w-{Y}^0_-(v,x_2)Y^0_-(u,x_1)w\nonumber\\
&=&Sing (Y^0_-(Y^0_-(u,x_1-x_2)v,x_2)w).
\end{eqnarray}
Let $\{a_i|i\in I\}$ and $\{b_j|j\in J\}$ be ordered bases of
$\mathbb{A}$ and $\mathbb{B}$, respectively. We assume that
$a_i\leq b_j$ for all $i\in I, j\in J$. Hence ${C}=\{a_i,b_j|i\in
I, j\in J\}$ is an ordered basis of $\mathbb{A}\oplus \mathbb{B}$.
Let $u\leq v\leq w\in {C}$. We claim that
\begin{eqnarray}
&&{Y}^0_-(u,x_1)Y^0_-(v,x_2)w-{Y}^0_-(v,x_2)Y^0_-(u,x_1)w\nonumber\\
&=&Sing (e^{x_2\D}{Y}^0_-(w,-x_2)Y^0_-(u,x_1-x_2)v).\label{vpa}
\end{eqnarray}

Case I: If $u,v,w\in \mathbb{A}$, by (\ref{yaa'}) all the three
terms are straightly zero.

Case II: If $u,v\in \mathbb{A}$, $w\in \mathbb{B}$, by
(\ref{yaa'}), (\ref{yab}) all the three terms are zero.

Case III: Assume that $u\in \mathbb{A}$, $v,w\in \mathbb{B}$. By
(\ref{yab})-(\ref{ybb'}), and (\ref{1tcs})-(\ref{1tca}), we have
\begin{eqnarray*}
&&{Y}^0_-(u,x_1)Y^0_-(v,x_2)w-{Y}^0_-(v,x_2)Y^0_-(u,x_1)w\\
&=&Y^0_-(u,x_1)(v_0wx_2^{-1}+v_1wx_2^{-2})-Y^0_-(v,x_2)(u_0wx_1^{-1})\\
&=&u_0v_0wx_1^{-1}x_2^{-1}-v_0u_0wx_1^{-1}x_2^{-1}\\
&=&(u_0v)_0wx_1^{-1}x_2^{-1},
\end{eqnarray*}
and
\begin{eqnarray*}
& &Sing
(e^{x_2\D}{Y}^0_-(w,-x_2)Y^0(u,x_1-x_2)v)\\
&=&Sing(e^{x_2\D}w_0u_0v(-x_2)^{-1}(x_1-x_2)^{-1})\\
&=&-w_0u_0vx_1^{-1}x_2^{-1}\\
&=&(u_0v)_0wx_1^{-1}x_2^{-1}.
\end{eqnarray*}
Hence, (\ref{vpa}) holds when $u\in \mathbb{A}, v,w\in\mathbb{B}$.

Case IV: Assume that $u,v,w\in \mathbb{B}$. By
(\ref{1tcd})-(\ref{1tca}), (\ref{yba})-(\ref{ybb'}), and
Proposition \ref{Dpartial}, we have
\begin{eqnarray*}
& &{Y}^0_-(u,x_1)Y^0_-(v,x_2)w-{Y}^0_-(v,x_2)Y^0_-(u,x_1)w\\
&=&Y^0_-(u,x_1)(v_0wx_2^{-1}+v_1wx_2^{-2})-Y^0_-(v,x_2)(u_0wx_1^{-1}+u_1wx_1^{-2})\\
&=&u_0v_0wx_1^{-1}x_2^{-1}+u_1v_0wx_1^{-2}x_2^{-1}+u_0v_1wx_1^{-1}x_2^{-2}\\
&
&-v_0u_0wx_1^{-1}x_2^{-1}-v_1u_0wx_1^{-1}x_2^{-2}-v_0u_1wx_1^{-2}x_2^{-1}\\
&=&(u_0v)_0wx_1^{-1}x_2^{-1}-(v_0u)_1wx_1^{-2}x_2^{-1}+(u_0v)_1wx_1^{-1}x_2^{-2},
\end{eqnarray*}
and
\begin{eqnarray*}
& &Sing
(e^{x_2\D}{Y}^0_-(w,-x_2)Y^0_-(u,x_1-x_2)v)\\
&=&Sing(e^{x_2\D}(w_0u_0v(-x_2)^{-1}(x_1-x_2)^{-1}+w_1u_0v(-x_2)^{-2}(x_1-x_2)^{-1}))\\
& &+Sing(e^{x_2\D}w_0u_1v(-x_2)^{-1}(x_1-x_2)^{-2})\\
&=&-w_0u_0vx_1^{-1}x_2^{-1}+w_1u_0vx_1^{-1}x_2^{-2}+w_1u_0vx_1^{-2}x_2^{-1}+\partial(w_1u_0v)x_1^{-1}x_2^{-1}\\
& &-w_0u_1vx_1^{-2}x_2^{-1}\\
&=&(u_0v)_0wx_1^{-1}x_2^{-1}+(u_0v)_1wx_1^{-1}x_2^{-2}+((u_0v)_1w-w_0u_1v)x_1^{-2}x_2^{-1}\\
&=&(u_0v)_0wx_1^{-1}x_2^{-1}+(u_0v)_1wx_1^{-1}x_2^{-2}+(-(v_0u)_1w+(\partial(v_1u))_1w+(u_1v)_0w)x_1^{-2}x_2^{-1}\\
&=&(u_0v)_0wx_1^{-1}x_2^{-1}+(u_0v)_1wx_1^{-1}x_2^{-2}+(-(v_0u)_1w-(v_1u)_0w+(u_1v)_0w)x_1^{-2}x_2^{-1}\\
&=&(u_0v)_0wx_1^{-1}x_2^{-1}+(u_0v)_1wx_1^{-1}x_2^{-2}-(v_0u)_1wx_1^{-2}x_2^{-1}.
\end{eqnarray*}
Hence, (\ref{vpa}) holds when  $u,v,w\in\mathbb{B}$.

Therefore, for $u,v,w\in C$ such that $u\leq v\leq w$, we have
$${Y}^0_-(u,x_1)Y^0_-(v,x_2)w-{Y}^0_-(v,x_2)Y^0_-(u,x_1)w=Sing
(e^{x_2\D}{Y}^0_-(w,-x_2)Y^0_-(u,x_1-x_2)v).$$ By Corollary
\ref{lhs} and Proposition \ref{singsing} ,we also have
\begin{eqnarray}
[Y^0_-(u,x_1),Y^0_-(v,x_2)]w&=&
Sing(e^{x_2\D}Sing(e^{-x_2\D}Y^0_-(Y^0_-(u,x_1-x_2)v,x_2)w)\nonumber\\
&=&Sing(Y^0_-(Y^0_-(u,x_1-x_2)v,x_2)w) \label{yuvw}\end{eqnarray}
for all $u,v,w\in C$ such that $u\leq v\leq w$. Moreover,
(\ref{yuvw}) holds for all $u,v\in\mathbb{A}\oplus \mathbb{B}$.

Next, we will show that for $u,v\in\mathbb{A}\oplus\mathbb{B}$,
$m,n\in\N$,
$$[Y^0_-(\D^nu,x_1),Y^0_-(\D^m v,x_2)]=Sing(Y^0_-(Y^0_-(\D^n
u,x_1-x_2)\D^m v,x_2)).$$ We will follow the last part of the
proof of Theorem 3.6 in \cite{li2}. Let $u,v\in \mathbb{A}\oplus
\mathbb{B}$. By Corollary \ref{exd}, we have
\begin{eqnarray*}
& &[Y^0_-(e^{z_1\D}u,x_1),Y^0_-(e^{z_2\D}v,x_2)]e^{z\D}w\\
&=&e^{z_1\frac{\partial}{\partial x_1}}e^{z_2\frac{\partial}{\partial x_2}}[Y^0_-(u,x_1),Y^0_-(v,x_2)]e^{z\D}w\\
&=&e^{(z_1-z)\frac{\partial}{\partial
x_1}}e^{(z_2-z)\frac{\partial}{\partial x_2}}e^{z\D}[Y^0_-(u,x_1),Y^0_-(v,x_2)]w\\
&=&e^{(z_1-z)\frac{\partial}{\partial
x_1}}e^{(z_2-z)\frac{\partial}{\partial
x_2}}e^{z\D}Sing(Y^0_-(Y^0_-(u,x_1-x_2)v,x_2)w)\\
&=&e^{(z_1-z)\frac{\partial}{\partial
x_1}}e^{(z_2-z)\frac{\partial}{\partial
x_2}}Sing_{x_1,x_2}(Y^0_-(Y^0_-(u,x_1-x_2)v,x_2+z)e^{z\D}w)\\
&=&Sing_{x_1,x_2}(Y^0_-(Y^0_-(u,x_1-x_2+z_1-z_2)v,x_2+z_2)e^{z\D}w)\\
&=&Sing_{x_1,x_2}(Y^0_-(Y^0_-(e^{(z_1-z_2)\D}u,x_1-x_2)v,x_2+z_2)e^{z\D}w)\\
&=&Sing_{x_1,x_2}(Y^0_-(e^{-z_2\D}Y^0_-(e^{z_1\D}u,x_1-x_2)e^{z_2\D}v,x_2+z_2)e^{z\D}w)\\
&=&Sing_{x_1,x_2}Y^0_-(Y^0_-(e^{z_1\D}u,x_1-x_2)e^{z_2\D}v,x_2)e^{z\D}w).
\end{eqnarray*}
Hence, $[Y^0_-(\D ^n u,x_1), Y^0_-(\D ^m v,x_2)]=Sing
(Y^0_-(Y^0_-(\D^n u,x_1-x_2)\D ^m v,x_2).$
\end{proof}

\begin{thm}\label{cvertexlie} $(\mathcal{C}, Y^0_-,\D)$ is a vertex Lie algebra. Moreover,
$S(\mathcal{C})$ is a vertex Poisson algebra. In particular,
$Y^0_-$ extend to a vertex Poisson algebra structure $Y_-$ on
$S(\mathcal{C})$ in the following way. First for
$u\in\mathcal{C}$, we define a unique element
$$\widetilde{Y^0_-}(u,x)\in x^{-1}(\Der S(\mathcal{C}))[[x^{-1}]]$$
by $$\widetilde{Y^0_-}(u,x)v=Y^0_-(u,x)v\ \ \ \text{ for }
u,v\in\mathcal{C}.$$ For $a\in S(\mathcal{C})$, we define
$$Y_-(a,x)\in x^{-1}(\Der S(\mathcal{C}))[[x^{-1}]]$$ by
$$Y_-(a,x)u=Sing(e^{x\D}\widetilde{Y^0_-}(u,-x)a)\ \ \ \text{ for }
u\in\mathcal{C}.$$
\end{thm}
\begin{proof} It follows from Proposition \ref{bracdy}, Corollary \ref{lhs}, Proposition \ref{bskew},
Proposition \ref{rvertexlie}, and Proposition 3.7 in \cite{li2}.
\end{proof}
\begin{coro}\label{grading2} Let $u\in \mathcal{C}$. We have $Y_-(u,x)v=Y^0_-(u,x)v $ for all $v\in\mathcal{C}$.
Furthermore, $Y_-(u,x)a=\widetilde{Y^0_-}(u,x)a$ for all $a\in
S(\mathcal{C})$.
\end{coro}

Since $\mathcal{C}$ is an $\N$-graded vector space, it implies
that $S(\mathcal{C})=\coprod_{n\in\N}S(\mathcal{C})_{(n)}$ is an
$\N$-graded vector space. By Corollary \ref{dui}, Theorem
\ref{cvertexlie} and Corollary \ref{grading2}, we can conclude
that:
\begin{lem}\label{grading3} For $u\in\mathcal{C}_{(m)}$, $i\in\N$,
$u_iS(\mathcal{C})_{(n)}\subset S(\mathcal{C})_{(n+m-i-1)}$. Here
$u_i\in \Der S(\mathcal{C})$.
\end{lem}

\begin{thm}\label{scvpa} $(S(\mathcal{C}),\D )$ is an $\N$-graded vertex
 Poisson
 algebra. In particular, we have $S(\mathcal{C})_{(0)}$ is $S(\mathbb{A})$
and for $n\geq 1$,
\begin{eqnarray*}
& &S(\mathcal{C})_{(n)}\\
&=&span_{\C}\{\ \ \D^{n_1}(b_1)\cdot
\D^{n_2}(b_2)\cdot....\cdot\D^{n_k}(b_k)\cdot a_1\cdot...\cdot
a_l\ \ |\ \ a_j\in\mathbb{A},\ \ b_i\in \mathbb{B},\\
& &\hspace{2cm}l\in\N,\ \ k\geq 1,\ \ n_1\geq...\geq n_k\geq 0,\ \
n_1+...+n_k+k=n\ \ \}.
\end{eqnarray*}
Here, $S(\mathbb{A})$ is the symmetric algebra over the space
$\mathbb{A}$.
\end{thm}

\begin{proof} First, we show that $\D S(\mathcal{C})_{(n)}\subset S(\mathcal{C})_{(n+1)}$ for $n\in \N$.
Clearly, for $n\in\N$, $\D
\mathcal{C}_{(n)}\subset\mathcal{C}_{(n+1)}$. Since $\D$ is a
derivation on $S(\mathcal{C})$, we can conclude immediately that
$\D S(\mathcal{C})_{(n)}\subset S(\mathcal{C})_{(n+1)}$.

Next, let $a\in S(\mathcal{C})_{(n)}$, $i\in\N$. We will show that
$$a_iS(\mathcal{C})_{(r)}\subset S(\mathcal{C})_{(r+n-i-1)}\ \
\text{ for all }r\in\N. $$ Let $r\in\N$ and let
$u\in\mathcal{C}_{(r)}$. By Corollary \ref{grading2} and Lemma
\ref{grading3}, we have
\begin{eqnarray*} &
&Y_-(a,x)u\\
&=&Sing(e^{x\D}\widetilde{Y^0}_-(u,-x)a)\\
&=&Sing(e^{x\D}\sum_{j=0}^{\deg a+\deg  u-1}u_ja(-x)^{-j-1})\\
&=&\sum_{j= 0}^{\deg a+\deg u-1}\sum_{k=0}^{j}(-1)^{-j-1}\frac{\D
^{j-k}}{(j-k)!}(u_{j}a)x^{-1-k}.
\end{eqnarray*}
Here, $u_j\in \Der S(\mathcal{C})$. It implies that
$a_i\mathcal{C}_{(r)}\subset {S(\mathcal{C})}_{(r+n-i-1)}$. Since
$a_i$ is a derivation on $S(\mathcal{C})$ and
$a_i\mathcal{C}_{(m)}\subset S(\mathcal{C})_{(m+n-i-1)}$ for all
$m\in\N$, it follows that $a_iS(\mathcal{C})_{(m)}\subset
S(\mathcal{C})_{(m+n-i-1)}$ for all $m\in\N$. Therefore,
$(S(\mathcal{C}),Y_-)$ is an $\N$-graded vertex Poisson algebra.

The second statement is clear.
\end{proof}

\section{Vertex Poisson algebras associated with Courant
algebroids}

In this section, we construct an $\N$-graded vertex Poisson
algebra
$S(\mathcal{C})_{\mathcal{B}}=\coprod_{n\in\N}(S(\mathcal{C})_{\mathcal{B}})_{(n)}$
associated with a Courant $\mathcal{A}$-algebroid $\mathcal{B}$.
Also, we show that $(S(\mathcal{C})_{\mathcal{B}})_{(0)}$ can be
naturally identified with $\mathcal{A}$ as a commutative
associative algebra and $(S(\mathcal{C})_{\mathcal{B}})_{(1)}$ can
be identified with $\mathcal{B}$ as a Courant
$\mathcal{A}$-algebroid.

For the rest of this paper, {\em we assume that $\mathcal{A}$ is a
commutative associative algebra with the identity $e$ and
$\mathcal{B}$ is a Courant $\mathcal{A}$-algebroid}.  By
Proposition \ref{cor1trun}, $\mathcal{A}\oplus \mathcal{B}$ is a
1-truncated conformal algebra such that for $a,a'\in \mathcal{A}$,
$u,v\in \mathcal{B}$,
\begin{eqnarray}
& &(au)_0a'=a(u_0a'),\label{Dera1}\\
& &(au)_1v=a(u_1v)=u_1(av),\label{Syma}\\
& &u_0(av)=a(u_0v)+(u_0a)v,\label{Dera2}\\
& &u_0(aa')=a(u_0a')+(u_0a)a'.\label{Dec}
\end{eqnarray}
Furthermore, by Theorem \ref{scvpa}, we have an $\N$-graded vertex
Poisson algebra
$S(\mathcal{C})=\coprod_{n\in\N}S(\mathcal{C})_{(n)}$ associated
with a 1-truncated conformal algebra $\mathcal{A}\oplus
\mathcal{B}$. We denote a multiplication on $S(\mathcal{C})$ by
"$\cdot$". Let
\begin{eqnarray}
E_0&=&span\{\ \ e-1,\ \ a\cdot a'-aa'\ \ |\ \ a,\ \
a'\in\mathcal{A}\ \ \}\subset
S(\mathcal{C})_{(0)},\\
E_1&=&span\{\ \ a\cdot b-ab \ \ |\ \ a\in\mathcal{A},\ \ b\in\mathcal{B}\ \ \}\subset S(\mathcal{C})_{(1)},\\
E&=&E_0\oplus E_1\subset S(\mathcal{C}).
\end{eqnarray}

\begin{lem}\label{pid} For $u\in \mathcal{A}\oplus \mathcal{B}$, $n\in\N$, we
have $u_nE\subset E$. Moreover, $\D E_0\subset E_1$.
\end{lem}
\begin{proof} We first show that for $u\in \mathcal{A}\oplus \mathcal{B}$, $n\in\N$, $$u_nE_0\subset E_0.$$
For $a\in\mathcal{A}$, $n\in\N$, we have $a_nE_0=0$. Similarly,
for $b\in\mathcal{B}$, $n\geq 1$, we have $b_nE_0=0$. Let
$b\in\mathcal{B}$. By Corollary \ref{dera}, Corollary \ref{dee},
and (\ref{Dec}), we have
$$b_0(e-1)=0$$ and
\begin{eqnarray*}
& &b_0(a\cdot a'-aa')\\
&=&a\cdot b_0(a')+b_0(a)\cdot a'-(b_0a)a'-a(b_0a')\\
&=&(a\cdot( b_0a')-a(b_0a'))+((b_0a)\cdot a'-(b_0a)a')\in E_0.
\end{eqnarray*}
for all $a,a'\in\mathcal{A}$. Therefore, for
$u\in\mathcal{A}\oplus\mathcal{B}$, $n\in\N$, $u_nE_0\subset E_0$.

Next, we show that for $u\in \mathcal{A}\oplus \mathcal{B}$,
$n\in\N$, $u_nE_1\subset E$. Clearly, for $a\in \mathcal{A}$,
$b\in\mathcal{B}$, $m\geq 1$, $n\geq 2$, we have $$a_m E_1\subset
E_1, \text{ and } b_n E_1\subset E_1.$$ Let $a', a\in
\mathcal{A}$, $b\in \mathcal{B}$. By (\ref{1tcs}), (\ref{Dera1}),
we have
\begin{eqnarray*}
& &a'_0(a\cdot b-ab)\\
&=&a'_0(a)\cdot b+a\cdot a'_0(b)+(ab)_0a'\\
&=&a\cdot a'_0(b)+a(b_0a')\\
&=&a\cdot (a'_0b)-a(a'_0b)\in E_0.
\end{eqnarray*}
Hence, $$a_nE_1\subset E\ \ \ \ \text{ for all }\ \ a\in
\mathcal{A},\ \ n\in\N.$$ Let $a\in\mathcal{A}$, $u,b\in
\mathcal{B}$. By (\ref{Syma}), (\ref{Dera2}), we have
$$u_0(a\cdot b-ab)=(u_0a)\cdot b+a\cdot (u_0b)-a(u_0b)-(u_0a)b\in
E_1,$$ and
$$u_1(a\cdot b-ab)=(u_1a)\cdot b+a\cdot (u_1b)-u_1(ab)=a\cdot u_1(b)-a(u_1b)\in
E_0.$$ It follows that for $b\in\mathcal{B}$, $n\in\N$,
$b_nE_1\subset E$.

 Next, we show that $\D E_0\subset E_1$. By Proposition \ref{Dpartial}, we have that $$\D(e-1)=\D
e=\partial(e)=0\in E_1.$$ Let $a,a'\in A$. By Propositon
\ref{Dpartial}, and the fact that
 $\D$ is a derivation on $S(\mathcal{C})$ and $\partial$ is a
 derivation from $\mathcal{A}$ to $\mathcal{B}$, we have
\begin{eqnarray*}
& &\D(a\cdot a'-aa')\\
&=&a\cdot (\D a')+(\D a)\cdot
a'-\D (aa')\\
&=&a\cdot \partial (a')+\partial( a)\cdot
a'-\partial(aa')\\
&=&a\cdot \partial (a')+a'\cdot\partial (a)
-a\partial(a')-a'\partial(a)\in  E_1.
\end{eqnarray*}
Therefore, $\D E_0\subset E_1$.
\end{proof}
For a subset $U$ of $S(\mathcal{C})$, we set $$\C[\D](U)=span\{\ \
\D ^m(u)\ \ |\ \ u\in U,\ \ m\in \N\ \ \}.$$ Define
$$I_{\mathcal{B}}= S(\mathcal{C})\cdot \C[\D](E)\subset S(\mathcal{C}),$$ an ideal of a
commutative associative algebra $S(\mathcal{C})$.
\begin{lem}\label{gradideal} $I_{\mathcal{B}}$ is an ideal of a vertex Poisson algebra
$S(\mathcal{C})$. Moreover,
\begin{eqnarray*}
& &S(\mathcal{C})_{(0)}=(I_{\mathcal{B}})_{(0)}\oplus
\mathcal{A}\text{ and,}\\
& &S(\mathcal{C})_{(1)}=(I_{\mathcal{B}})_{(1)}\oplus
\mathcal{B}.\end{eqnarray*}
\end{lem}
\begin{proof} Since $\D$ is a derivation on $S(\mathcal{C})$ and
$\D(\C[\D](E))\subset\C[\D](E)$, it follows that
\begin{equation}\label{dib}\D I_{\mathcal{B}}\subset
I_{\mathcal{B}}.\end{equation} Let
$v\in\mathcal{A}\oplus\mathcal{B}$, and $u\in E$. We claim that
for $m\in\N$,
$$v_i\D ^m(u)\in\C[\D](E)\ \ \text{ for all }i\in\N.$$ We prove this by an induction on $m$. For
the case when $m=0$, it follows immediately from Lemma \ref{pid}.
Recall that for $i\in\N$, $[\D, v_i]=-iv_{i-1}$ (see Proposition
\ref{rhd}). We now assume that for $n< m$, $v_i\D
^n(u)\in\C[\D](E)$ for all $i\in\N$. Observe that $$v_i\D ^m(u)=\D
v_i\D^{m-1}(u)+i v_{i-1}\D^{m-1}(u).$$ By an induction hypothesis,
we conclude that $v_i\D^m(u)\in \C[\D](E)$ for all $i\in\N$.
Therefore, for $m\in\N$,
\begin{equation}\label{vid}v_i\D^m(u)\in\C[\D](E)\text{ for all } i\in\N.\end{equation}

Next, we show that for $v\in \mathcal{A}\oplus\mathcal{B}$,
$i\in\N$, $v_iI_{\mathcal{B}}\subset I_{\mathcal{B}}$. Let
$v\in\mathcal{A}\oplus\mathcal{B}$ and $i\in\N$. By (\ref{hd}),
(\ref{vid}), we have that for $a\in S(\mathcal{C})$, $u\in E$,
$m\in\N$,
$$v_i(a\cdot \D^m u)=a\cdot v_i(\D^m u)+v_i(a)\cdot \D ^m u\in
I_{\mathcal{B}}.$$ This implies that for $v\in
\mathcal{A}\oplus\mathcal{B}$, $i\in\N$,
\begin{equation}\label{viib} v_iI_{\mathcal{B}}\subset
I_{\mathcal{B}}.\end{equation}

Next, we will show that for $a\in S(\mathcal{C})$, $w\in
I_{\mathcal{B}}$, $i\in\N$, $$a_iw\in I_{\mathcal{B}}.$$ By
(\ref{hs}), (\ref{dib}), (\ref{viib}), we have that for $w\in
I_{\mathcal{B}}$, $v\in\mathcal{A}\oplus\mathcal{B}$, $i\in\N$,
$$w_iv\in I_{\mathcal{B}}.$$ Furthermore, by (\ref{viib}),
Corollary \ref{dui}, (\ref{hs}), and (\ref{dib}), we conclude that
$$(\D^mv)_i w\in I_{\mathcal{B}},\ \ \text{ and } w_i(\D^mv)\in
I_{\mathcal{B}}$$ for all $v\in\mathcal{A}\oplus\mathcal{B}$,
$w\in I_{\mathcal{B}}$, $i\in \N$, and $m\geq 1$. Hence, for $w\in
I_{\mathcal{B}}$, $i\in\N$, $$w_i\mathcal{C}\subset
I_{\mathcal{B}}.$$ Let $w\in I_{\mathcal{B}}$, $i\in \N$. Since
$w_i$ is a derivation on $S(\mathcal{C})$ and
$w_i\mathcal{C}\subset I_{\mathcal{B}}$, it follows that
\begin{equation}\label{wisc}w_iS(\mathcal{C})\subset I_{\mathcal{B}}.\end{equation} Hence, by
(\ref{hs}), (\ref{dib}), (\ref{wisc}), we can conclude that for
$a\in S(\mathcal{C})$, $w\in I_{\mathcal{B}}$, $ i\in\N$, $a_i
w\in I_{\mathcal{B}}$ and $I_{\mathcal{B}}$ is an ideal of
$S(\mathcal{C})$ as a vertex Poisson algebra.

Next, we will show that
$S(\mathcal{C})_{(0)}=(I_{\mathcal{B}})_{(0)}\oplus \mathcal{A}$
and $S(\mathcal{C})_{(1)}=(I_{\mathcal{B}})_{(1)}\oplus
\mathcal{B}.$ Clearly,
$$(I_{\mathcal{B}})_{(0)}=S(\mathcal{A})E_0\text{ and
}(I_{\mathcal{B}})_{(1)}=S(\mathcal{A})\D E_0+S(\mathcal{A})\cdot
\mathcal{B}\cdot E_0+S(\mathcal{A})\cdot E_1.$$ By Lemma
\ref{pid}, we have $S(\mathcal{A})\D E_0\subset
S(\mathcal{A})\cdot E_1.$ Since $\mathcal{B}$ is an
$\mathcal{A}$-module, it implies that for $a,a'\in\mathcal{A}$,
$b\in\mathcal{B}$
\begin{eqnarray*} & &b\cdot (a\cdot
a')-b\cdot
(aa')\\&=&(a\cdot a')\cdot b-(aa')\cdot b\\
&=&a\cdot (a'\cdot b)-a\cdot (a'b)+a\cdot (a'b)-(aa')\cdot b\\
&=&a\cdot (a'\cdot b)-a\cdot (a'b)+a\cdot (a'b)-(aa')\cdot
b+(aa')b-(aa')b\\
&=&a\cdot (a'\cdot b)-a\cdot (a'b)+a\cdot (a'b)-(aa')\cdot
b+(aa')b-a(a'b)\in S(\mathcal{A})\cdot E_1. \end{eqnarray*}
Therefore, $$S(\mathcal{A})\cdot \mathcal{B}\cdot E_0\subset
S(\mathcal{A})\cdot E_1, \text{ and
}(I_{\mathcal{B}})_{(1)}=S(\mathcal{A})\cdot E_1.$$ Moreover, we
have
$$S(\mathcal{C})_{(0)}=S(\mathcal{A})=\mathcal{A}\oplus
(I_{\mathcal{B}})_{(0)},$$ and
$$S(\mathcal{C})_{(1)}= S(\mathcal{A})\cdot \mathcal{B}=\mathcal{B}\oplus(I_{\mathcal{B}})_{(1)}.$$
\end{proof}

Set
$$S(\mathcal{C})_{\mathcal{B}}=S(\mathcal{C})/I_{\mathcal{B}},$$ an
$\N$-graded vertex Poisson algebra.
\begin{thm} Let $\mathcal{B}$ be a Courant
$\mathcal{A}$-algebroid and let
$S(\mathcal{C})_{\mathcal{B}}=\coprod_{n\in\N}\left(S(\mathcal{C})_{\mathcal{B}}\right)_{(n)}$
be the associated $\N$-graded vertex Poisson algebra. We have that
$\left(S(\mathcal{C})_{\mathcal{B}}\right)_{(0)}=\mathcal{A}$ and
$\left(S(\mathcal{C})_{\mathcal{B}}\right)_{(1)}=\mathcal{B}$.
Moreover, for any $n\geq 1$,
\begin{eqnarray*}
\left(S(\mathcal{C}_{\mathcal{B}})\right)_{(n)}
&=&span\{\ \ \D^{n_1}(b_1)\cdot\D^{n_2}(b_2)\cdot...\cdot\D^{n_k}(b_k)\\
& &\hspace{2cm} |\ \ b_i\in\mathcal{B}, \ \ n_1\geq...\geq n_k\geq
0,\ \ k\geq 1,\ \ n_1+...+n_k+k=n\ \ \}.
\end{eqnarray*}
\end{thm}
\begin{proof} By Lemma \ref{gradideal}, we have $\left(S(\mathcal{C})_{\mathcal{B}}\right)_{(0)}=\mathcal{A}$ and
$\left(S(\mathcal{C})_{\mathcal{B}}\right)_{(1)}=\mathcal{B}$. For
a subset $U$ of $S(\mathcal{C})_{\mathcal{B}}$, we set
$$\C[\D]U=Span\{\ \ \D^m(u)\ \ |\ \ u\in U,\ \
m\in\N\ \ \}.$$ We will show that
$$\coprod_{n\geq
1}\left(S(\mathcal{C})_{\mathcal{B}}\right)_{(n)}=S(\C[\D]\mathcal{B})\cdot
\C[\D]\mathcal{B}.$$ Here, $S(\C[\D]\mathcal{B})$ is a symmetric
algebra over the space $\C[\D]\mathcal{B}$. By Theorem
\ref{scvpa}, it is enough to show that
$(\C[\D]\mathcal{B})\cdot\mathcal{A}\subset\C[\D]\mathcal{B}.$
Clearly, for $a\in\mathcal{A}$, $b\in \mathcal{B}$, we have
$$b\cdot a=a\cdot b=ab+a\cdot b-ab=ab\in \mathcal{B}.$$ Recall that for
$u,v\in S(\mathcal{C})_{\mathcal{B}}$, $n\in\N$, we have
$$\D^n(u\cdot v)=\sum_{i=0}^n{n\choose i}\D^i(u)\cdot \D^{n-i}(v).$$
Let $a\in\mathcal{A}$, $b\in \mathcal{B}$, $n\geq 1$. By
Proposition \ref{Dpartial}, we have
\begin{eqnarray*}
(\D^n b)\cdot a &=&\D^n(a\cdot b)-\sum_{i=1}^n{n\choose
i}\D^i(a)\cdot
\D^{n-i}(b)\\
&=&\D^n(ab)-\sum_{i=1}^n{n\choose i}\D^{i-1}(\partial (a))\cdot
\D^{n-i}(b).
\end{eqnarray*}
Therefore, for $a\in\mathcal{A}$, $b\in\mathcal{B}$ and $n\in\N$,
we have $(\D^n b)\cdot a\in S(\C[\D]\mathcal{B})\cdot
\C[\D]\mathcal{B}$. It follows that$$\coprod_{n\geq
1}\left(S(\mathcal{C})_{\mathcal{B}}\right)_{(n)}=S(\C[\D]\mathcal{B})\cdot
\C[\D]\mathcal{B}.$$
\end{proof}

\section{Appendix}
\begin{de}\cite{lli} A {\em vertex algebra} is a vector space $V$ equipped with a
linear map
\begin{eqnarray}
Y: V&\rightarrow& (\End V)[[x,x^{-1}]],\nonumber\\
v&\mapsto& Y(v,x)=\sum_{n\in\Z}v_{n}x^{-n-1} \;\;\;(\mbox{where
}v_{n}\in \End V)
\end{eqnarray}
and equipped with a distinguished vector ${\bf 1}$, called the
{\em vacuum (vector)}, such that for $u,v\in V$,
\begin{eqnarray}
& &u_{n}v=0\;\;\;\mbox{ for n sufficiently large},\\
& &Y({\bf 1},x)=1,\\
& &Y(v,x){\bf 1}\in V[[x]]\;\;\mbox{ and } \lim_{x\rightarrow
0}Y(v,x){\bf 1}=v
\end{eqnarray}
and
\begin{eqnarray}
&
&x_{0}^{-1}\delta\left(\frac{x_{1}-x_{2}}{x_{0}}\right)Y(u,x_{1})Y(v,x_{2})
-x_{0}^{-1}\delta\left(\frac{x_{2}-x_{1}}{-x_{0}}\right)Y(v,x_{2})Y(u,x_{1})\nonumber\\
& &\ \ \ \
=x_{2}^{-1}\delta\left(\frac{x_{1}-x_{0}}{x_{2}}\right)Y(Y(u,x_{0})v,x_{2})
\end{eqnarray}
the {\em Jacobi identity}.
\end{de}

\end{document}